\documentclass{amsart}
\theoremstyle{plain}

\newtheorem{thm}{Theorem}
\theoremstyle{definition}

\newtheorem{defn}[thm]{Definition}
\newtheorem{lem}[thm]{Lemma}
\newtheorem{rem}[thm]{Remark}
\newtheorem{prop}[thm]{Proposition}
\newtheorem{cor}[thm]{Corollary}

\usepackage[]{graphicx}
\usepackage[T1]{fontenc}
\usepackage[latin9]{inputenc}
\usepackage{color}
\usepackage{amstext}
\usepackage{amsthm}
\usepackage{amssymb}

\numberwithin{equation}{section}
\numberwithin{figure}{section}

\makeatletter

\begin{document}

\title{A note on Andrews inequalities}

\author{Fang Hao}
\address{The University of Iowa, Department of Mathematics, 
Iowa City, IA 52240, USA}
\email{hao-fang@uiowa.edu}

\author{Biao Ma}%
\email{biaoma@uiowa.edu}

\author{Wei Wei}%
 \address{Shanghai Center for Mathematical Sciences, 2005 Songhu Road, Shanghai 200438, China}

\email{wei\_wei@fudan.edu.cn }

\begin{abstract}
In this note, we prove two Andrews type inequalities concerning conic
manifolds and manifolds with boundary. In each case, we establish
a rigidity result when the equality holds.
\end{abstract}

\maketitle

\section{Introduction}

In his unpublished work, Ben Andrews has proved the following geometric
inequality for closed Riemannian manifolds with positive Ricci curvature:
\begin{thm}
\label{thm:original}Let $(M,g)$ be an $n$-dimensional closed Riemannian
manifold with positive Ricci curvature, $Ric_{g}.$ Then for any smooth
function $\phi\in C^{\infty}(M,g)$ such that $\int_{M}\phi dv_{g}=0$
\begin{equation}
\frac{n}{n-1}\int_{M}\phi^{2}\ dv_{g}\leq\int_{M}Ric^{-1}(\nabla\phi,\nabla\phi)\ dv_{g},\label{eq:original}
\end{equation}
where the equality holds for a non-vanishing $\phi$ if and only if
$(M,g)$ is conformal to the $n$-sphere $(S^{n},g_{S^{n}})$ with
its standard metric. In addition, $g$ is invariant under an $O(n-1)$
isometry group. 
\end{thm}

Inequality (\ref{eq:original}) has been proved and seen applications
in the study of Ricci flow, see section 3, Appendix B in \cite{BLN}.
More recently, Gursky and Streets \cite{Gursky-Streets1,Gursky-Streets2}
have used (\ref{eq:original}) in a different area of geometric analysis,
namely, conformal geometry of 4-manifolds. The $\sigma_{k}$-Yamabe
problem is a recent focus in the study of conformal geometry. It is
well known that in dimension 4, the positivity of $\sigma_{1}$ and
$\sigma_{2}$-curvatures implies the positivity of Ricci curvature.
Thus, (\ref{eq:original}) is applicable to such manifolds. Motivated
by a similar construction in K$\ddot{a}$hler geometry, Gursky-Streets
\cite{Gursky-Streets1} have defined a formal metric structure on
all conformal metrics over a 4-manifold $M$ with positive $\sigma_{2}$-curvature.
They then applied (\ref{eq:original}) to establish the geodesic convexity
of a crucial functional defined by Chang-Gursky-Yang \cite{CGY1}.
Through further careful analysis, Gursky and Streets are able to prove
the uniqueness of the $\sigma_{2}$-Yamabe problem solution. This
is a surprising result in contrast with the general dimension situation
for the standard Yamabe problem.

In this note, we prove some extensions of Andrews inequality (\ref{eq:original})
and explore further applications in conformal geometry. First we discuss
the case of conic manifolds. 
\begin{defn}
\label{def:conic manifolds}Let $(M^{n},g_{0})$ be a compact smooth
Riemannian $n$-manifold. For some $k\in\mathbb{N},$ assume that
$p_{i}\in M$ and $0>\beta_{i}>-1$ for $1\leq i\leq k$. Define a
conformal divisor 
\[
D=\sum_{i=1}^{k}p_{i}\beta_{i}.
\]
Let $M\backslash D=M\backslash\{p_{i}:i=1,2,\cdots,k\}$. Let $\gamma(x)$
be a function in $C^{\infty}(M^{n}\backslash D)$ such that near $p_{i}$,
$\gamma(x)-\beta_{i}\log r_{i}$ is locally smooth, where $r_{i}=\mathrm{dist}_{g_{0}}(x,p_{i})$
is the distance to $p_{i}$. Let 
\begin{equation}
g_{D}=e^{2\gamma}g_{0}.\label{eq:g_D definition}
\end{equation}
$g_{D}$ is a metric conformal to $g_{0}$ on $M^{n}\backslash D$
which has a conical singularity at each $p_{i}$. For a technical
reason, we define the following function spaces: 
\[
\tilde{C}_{\beta}^{2}(B_{R}(0),g_{0}):=\{w\in C^{\infty}(B_{R}(0)\backslash\{0\})\cap C^{0}(B_{R}(0)):|r^{-k\beta}\nabla^{k}w|<\infty,k=0,1,2\},
\]
and
\[
\tilde{C}_{D}^{2}(M,g_{0}):=\{w\in C^{\infty}(M\backslash D,g_{0}):\exists\delta>0,w|_{B_{\delta}(p_{i})}\in\tilde{C}_{\beta_{i}}^{2}(B_{\delta}(p_{i}),g_{0})\}.
\]
Then, we define the conformal class of $g_{D}$
\[
[g_{D}]:=\{g_{w}=e^{2w}g_{D}:w\in\tilde{C}_{D}^{2}(M,g_{0})\}.
\]
Note $[g_{D}]$ depends only on $(M^{n},g_{0})$ and $D$. We call
a $4$-tuple $(M^{n},g_{0},D,g_{1})$ a \textbf{conic $n$-manifold},
where\textbf{ $g_{0}$ }is the background metric and $g_{1}\in[g_{D}]$
is the conic\textbf{ }metric.
\end{defn}

We remark here that our definition of conic singularity is in the
sense of Cheeger-Colding \cite{Cheeger-Colding,Cheeger-Colding-1},
which is different from the similar notation used in studies of K$\ddot{a}$hler
geometry. For discussion of related conformal geometry problems on
conic manifolds, see for example \cite{FW1,FW2,FW3} and \cite{FM}.
We also note that for analytical consideration, we may relax the regularity
of the conformal factor function.

Our first main result is the following
\begin{thm}
\label{thm:(Conic-Andrews)} Let $(M,g_{0},D,g)$ be
a compact $n$-dimensional conic manifold with positive Ricci curvature.
Let $\phi\in C^{\alpha}(M,g_{0})\cap C^{1}(M\backslash D,g_{0})$
for some $\alpha>0$ such that $\int_{M}\phi dv_{g}=0$ . Then we
have
\begin{equation}
\frac{n}{n-1}\int_{M\backslash D}\phi^{2}\ dv_{g}\leq\int_{M\backslash D}Ric^{-1}(\nabla\phi,\nabla\phi)\ dv_{g}.\label{eq:conic}
\end{equation}
In particular, the equality holds for a non-vanishing $\phi$ if and
only if $(M,g_{0},D,g)$ is conformal to a conic sphere with at most $2$ conic points and an $O(n-1)$
isometry group, which is described in Theorem \ref{thm:warped product-1}.
\end{thm}

The second part of our note is on manifolds with boundary. Our result
is the following
\begin{thm}
\label{thm:andrews ineq with Boundary }Let $(M,g)$ be a connected
$n$-dimensional compact manifold with smooth boundary $\partial M$.
Suppose that $M$ has positive Ricci curvature. Let ${\rm II}$ be
the second fundamental form of $\partial M,$ and assume that ${\rm II}\geq0$
or equivalently, $\partial M$ is weakly convex. Then, for any $\phi$
such that $\int_{M}\phi dv_{g}=0$, 
\begin{equation}
\frac{n}{n-1}\int_{M}\phi^{2}\ dv_{g}\leq\int_{M}Ric^{-1}(\nabla\phi,\nabla\phi)\ dv_{g},\label{eq:with boundary}
\end{equation}
where the equality holds for a non-vanishing $\phi$ if and only if
$(M,g)$ is conformal to an n-dimensional hemisphere $S_{+}^{n}$
with its standard metric; furthermore, $g$ admits an $O(n-1)$ isometry
group. 
\end{thm}

\begin{rem}
In \cite{Escobar}, Escobar has studied the Laplacian eigenvalues
with Neumann boundary condition for manifolds with $Ric\geq(n-1)$
and convex boundary. He proved that the first eigenvalue $\lambda_{1}\geq n$
with the equality holding if and only if the manifold is isometric
to the standard hemisphere. In our case, $M$ is less restricted when
the equality in (\ref{eq:with boundary}) holds. 
\end{rem}

Theorem \ref{thm:(Conic-Andrews)} and Theorem \ref{thm:andrews ineq with Boundary }
are natural extensions of the original Andrews inequality, Theorem
\ref{thm:original}. We note that for various Andrews inequalities
we have obtained, the equality cases all imply the conformal-flat-ness
and the existence of a relatively large isometry group.

We make some further comments. First, we raise a general question
which may be of interest. For manifolds with positive Ricci curvature,
we may define a functional for non-zero square integrable functions
$\int_{M}\phi dv_{g}=0$:
\[
F(\phi)=\frac{\int_{M}Ric^{-1}(\nabla\phi,\nabla\phi)\ dv_{g}}{\int_{M}\phi^{2}\ dv_{g}}
\]
 and (\ref{eq:original}) obviously gives a lower bound estimate of $F$.
The corresponding Euler equation is thus
\begin{equation}
\delta(Ric^{-1}d\phi)=\lambda\phi,\label{eq:first eigenvalue ric inverse}
\end{equation}
which can be viewed as an eigenvalue problem. 

Obviously, when the manifold is Einstein, this is just the well-known
eigenvalue problem for the Laplacian operator. However, there are
known examples of non-Einstein manifolds with positive Ricci curvature.
For instance, if $k\geq4$, the connected sum of complex projective
spaces $M=\#_{k}\mathbb{CP}^{2}$ does not admit Einstein metric \cite{Besse},
while there always exists a metric on $M$ with positive Ricci curvature
due to Perelman \cite{Perelman}. Sha-Yang \cite{Sha-Yang,Sha-yang2}
has constructed a class of Ricci positive manifolds which admit arbitrarily
large Betti numbers. It may be of interest to study sharp Andrews
type inequality in these special cases.

Another direction of interest is to study the almost sharp Andrews
inequality. We may expect some rigidity results. For instance, when
the first eigenvalue in (\ref{eq:first eigenvalue ric inverse}) is
almost $\frac{n}{n-1}$, one may naively wish the manifold to be a
sphere. However, for dimension $n$ being even, Anderson \cite{Anderson} constructs
a metric on complex projective space $\mathbb{CP}^{n/2}$ with $Ric\geq n-1$
such that the Laplacian first eigenvalue is arbitrarily close to $n+1$.
We easily check that the first eigenvalue of (\ref{eq:first eigenvalue ric inverse})
in Anderson's example can be arbitrarily close to $\frac{n}{n-1}$.
For further discussions on eigenvalue problems and rigidity results,
see \cite{Cheng,Colding,petersen2}. 

The rest of the paper is organized as follows. In Section 2, we discuss
the warped product structure which is crucial in our study of sharp
Andrews inequality. In Section 3, we prove Theorem \ref{thm:(Conic-Andrews)}.
In Section 4, we prove Theorem \ref{thm:andrews ineq with Boundary }. 

\section{Warped product}

In this section, we discuss the geometry of warped products. Most
of the results in this section are well known and can be found in
\textcolor{black}{\cite{petersen,IshiharaTashiro,Tashiro,Cheeger-Colding-1,Wu-Ye,CMM,CSZ}.
We state them for future use.}
\begin{defn}
\label{def:Definition warped products}Let $(M^{n},g)$ be a $n$-dimensional
Riemannian manifold. We say $M$ is a \textbf{warped product} if $M\simeq[a,b]\times N^{n-1}$
and $g$ can be written as $g=dr^{2}+f(r)^{2}g_{N}$ for some $n-1$-dimensional
Riemannian manifold $(N^{n-1},g_{N})$ and a positive function $f:[a,b]\to\mathbb{R}$
. We use the following notation to denote a warped product:
\[
M=[a,b]\times_{f}N^{n-1}.
\]
\end{defn}

Let $\phi_{r}:N^{n-1}\hookrightarrow M^{n}$ be the embedding of $N$
at $r\in(a,b)$. Let $x^{i}$ be a local coordinate on $N^{n-1}$.
Then, $(r,x^{i})$ gives a coordinate on $M^{n}$. Let $g,\nabla,R_{ijkl},Ric$
be the metric, connection, curvature tensor and Ricci curvature tensor
on $M^{n}$ respectively and let $\bar{g}=g_{N},\bar{\nabla},\bar{R}_{ijkl},\bar{R}ic$
be the corresponding metric, connection, curvature tensor and Ricci
curvature tensor on $N$, respectively. By simple calculation, we
have:
\begin{lem}
\label{lem:petersen}
\begin{align*}
\nabla_{\partial_{r}}\partial_{i} & =(\ln f(r))'\partial_{i},\ \nabla_{\partial_{i}}\partial_{r}=(\ln f(r))'\partial_{i},\\
\nabla_{\partial_{i}}\partial_{j} & =-f(r)f'(r)\bar{g}_{ij}\partial_{r}+\bar{\nabla}_{\partial_{i}}\partial_{j}.
\end{align*}
The second fundamental form ${\rm II}=(h_{ij})$ of the embedding
$\phi_{r}$ is given by 
\[
{\rm II}=-f(r)^{-1}f'(r)\bar{g}_{ij}.
\]
By Gauss-Codazzi equation, we have
\[
R_{ijkl}-f(r)^{2}\bar{R}_{ijkl}=h_{ik}h_{jl}-h_{il}h_{jk},R_{rirj}=\left((\ln f)''+\left[(\ln f)'\right]^{2}\right)g_{ij}.
\]
For Ricci curvature tensor, we have
\end{lem}

\begin{align}
Ric_{ij} & =\left(-\frac{f''}{f}-(n-2)\left(\frac{f'}{f}\right)^{2}\right)g_{ij}+\overline{Ric}_{ij},\label{eq:Ricc}\\
Ric_{rr} & =-(n-1)\left(f''/f\right),Ric_{ir}=0.\nonumber 
\end{align}
The next proposition characterizes the warped product structure, see
\cite{petersen}.
\begin{prop}
\textcolor{black}{If there exists a function $u$ on $M$ such that
\begin{equation}
\nabla^{2}u-\frac{\Delta u}{n}g=0,\label{eq:gradient conformal vector field}
\end{equation}
and $du\not=0$, then $M$ is locally a warped product. If $du(p)=0$
at $p\in M$ and $\Delta u(p)\not=0$, then $g=dr^{2}+f(r)^{2}ds_{n-1}^{2}$,
where $ds_{n-1}^{2}=g_{S^{n-1}}$ is the standard metric on $(n-1)$-sphere.\label{prop:conformal vector field gives warped product}}
\end{prop}

\begin{proof}
We sketch the proof here. Let $s=\frac{du}{|du|}$, where $|du|$
is computed using $g^{T^{*}M}$, a metric on $T^{*}M$ that is canonically
induced from $g$. By (\ref{eq:gradient conformal vector field}),
we get $ds=0.$ Thus, we may define a local function $r$ such that
$dr=s.$ Thus, we have 
\begin{equation}
dr=\frac{du}{|du|}.\label{eq:add1}
\end{equation}
Let $H=\{x\in M,r=r_{0}\}$ be a level-set of $r$. Thus, for $x\in M$
near $H$, since $|dr|=1,$ the distance from $x$ to $H$ can be
computed as $|r-r_{0}|.$ By (\ref{eq:add1}), for any vector field
$V$ tangent to $TH$, we have $du(V)=0.$ Thus, $u$ is a function
of $r$. Let $f(r)=|\nabla u(r)|$. The metric $g$ on $M$ can then
be locally written as $g=dr^{2}+f(r)^{2}g_{H,r}$ in a open domain
$(r_{0}-\epsilon,r_{0}+\epsilon)\times H$ in $M$. Here $g_{H,0}$
is a metric on $H$. A direct computation shows that $g_{H,r}$ is
invariant along the $\partial_{r}$ direction. Hence $g_{H,r}=g_{H}$.
For a local critical point $p\in M$ of $u,$ it is easy to see that
$H$ has to be $S^{n-1}$. 
\end{proof}
\textcolor{black}{}
\begin{prop}
\textcolor{black}{If there exists a function $u$ on $M$ such that
\begin{equation}
\nabla^{2}u-\frac{\Delta u}{n}g=0,\label{eq:gradient conformal vector field-1}
\end{equation}
}and the Ricci curvature is positive, then $u$ is a Morse function,
i.e. critical points of $u$ are non-degenerate. Moreover, critical
points of $u$ are local maximums or minimums.\label{prop:critical points have trivial index}
\end{prop}

\begin{proof}
Since Ricci curvature is positive, by (\ref{eq:Ricc}), $f''/f<0$.
Hence, $f(r)$ is strictly concave. Therefore, along each geodesic
generated by $\partial_{r}$, $u(r)$ has at most one maximum and
one minimum. Recall $f(r(p))=|\nabla u(p)|$. If $p$ is a critical
point of $u,$ then $f(r(p))=0$ and
\[
\nabla^{2}u=f'(r)g,
\]
thus $p$ can be a local maximum, minimum, or a degenerate critical
point. Now, suppose that $p$ is a degenerate critical point of $u$.
Let $r_{1}=r(p)$. Then, $f(r_{1})=0$ and $f'(r_{1})=0$, hence $f(r)<0$
in a neighborhood of $r_{1}$ which contradicts with the fact that
$f$ is non-negative. Thus, $p$ is non-degenerate.
\end{proof}
Since the flow of $\nabla u$ gives a deformation retraction between
level-sets without crossing critical points, as an easy application
of Morse theory, we derive the following result.
\begin{cor}
\label{cor:warped sphere}Let $(M,g)$ be a connected compact manifold
with positive Ricci curvature. If there is a non-constant function
$u$ satisfying (\ref{eq:gradient conformal vector field}), then
$M$ is diffeomorphic to a sphere with warped product structure 
\[
g=dr^{2}+f(r)^{2}g_{S^{n-1}},
\]
where $g_{S^{n-1}}$ is the standard metric on $S^{n-1}$.
\end{cor}

\begin{rem}
By Morse theory, if $M$ only has maximum and minimum as its critical
points, then $M$ is homeomorphic to a sphere. The diffeomorphism
follows from the warped product structure.
\end{rem}

\textcolor{black}{A similar argument is used to derive the following
result for conic manifolds.}
\begin{thm}
\label{thm:warped product-1}Let $(M,g_{0},D,g)$ be a  conic manifold
with positive Ricci curvature where $D=\sum_{i=1}^{k}\beta_{i}p_{i}$.
 If there exists a non-constant function $u\in C^{1,\gamma}(B_{\delta}(p_{i}),g)\cap C^{\infty}(M\backslash D,g_{0})$
for some $\gamma>0$ such that 

\[
\nabla^{2}u=\frac{\Delta u}{n}g,\quad M\backslash D,
\]
then $|\nabla u|$ is constant on each level-set of $u.$ Let $f(r)=|\nabla u|$.
Then $(M,g)$ has to be one of the following: 

$A:$ $g=dr^{2}+f(r)^{2}g_{S^{n-1}}$ on $M\backslash\{q\}=S^{n-1}\times(0,\infty)$
with $\lim_{r\rightarrow0}\frac{f(r)}{r}=1+\beta$, here $\beta\in(-1,0)$
\textup{is} the conic coefficient and $D=\beta q$.

$B:$ $g=dr^{2}+f(r)^{2}g_{S^{n-1}}$ on $M\backslash\{q,p\}=S^{n-1}\times(a,b)$
with $\lim_{r\rightarrow b}|\frac{f(r)}{r-b}|=1+\beta_{b}$, here
$\beta_{b}\in(-1,0)$\textup{, $D=\beta_{b}p$} and $\lim_{r\rightarrow a}|\frac{f(r)}{r-a}|=1$\textup{. }

$C:$ $g=dr^{2}+f(r)^{2}g_{S^{n-1}}$ on $M\backslash\{q,p\}=S^{n-1}\times(a,b)$
with $\lim_{r\rightarrow b}|\frac{f(r)}{r-b}|=1+\beta_{b}$ and $\lim_{r\rightarrow a}|\frac{f(r)}{r-a}|=1+\beta_{a},$
\textup{here $\beta_{a}$, }$\beta_{b}\in(-1,0)$\textup{ and $D=\beta_{a}q+\beta_{b}p$.}
\end{thm}

\begin{proof}
We follow the proof of Proposition \ref{prop:conformal vector field gives warped product}
to define $r,r_{0},$ $H$ and $f(r)$. Hence, $H$ is a hypersurface
in $M$. We consider $\gamma(r)$, a maximally extended flow line
of $\frac{\partial}{\partial r}=\nabla u/|\nabla u|$. Thus, $\gamma(r)$
is a geodesic and $|\gamma'(r)|=1$. Since $M$ is compact, $\gamma$
is of finite length. Furthermore, $u|_{\gamma}$ is increasing with
respect to $r$. Also, recall
\[
f(r)=|\nabla u(r)|.
\]
Let $p_{1}$, $p_{2}\in\partial\{\gamma(r)\}$. Then, we have $f(r(p_{i}))=0,i=1,2.$
Therefore, $p_{i}$ is either a critical point of $u$, or a conic
point. We claim that there are finitely many choices of such $p_{i}$
for all possible choices of $\gamma$.

Fix $i\in\{1,2\}$. If $p_{i}$ is one of the conic points, then,
by definition they are isolated. We claim that, near $p_{i}$, the
metric $g$ can be written as $g=dr^{2}+f(r)^{2}g_{S^{n-1}}$, where
$f(r)\geq0$, $f(0)=0$. 

Near conic point $p_{i}$, let $g=e^{2v}d_{g_{0}}(x,p_{i})^{2\beta_{i}}g_{0}$,
where $v$ is in $\tilde{C}_{\beta_{i}}^{2}(B_{\delta}(p_{i}),g_{0})$
for some $\delta>0$. By Lemma \ref{lem:geodesic normal coordinate}
in the appendix, there is a geodesic normal coordinate $(\rho,\theta)$
such that $g=d\rho^{2}+h^{2}(\rho,\theta)g_{S^{n-1}}$ and $h=(1+\beta_{i})\rho+O(\rho^{2})$.
We have $\lim_{\rho\rightarrow0}h^{2}g_{S^{n-1}}/\rho^{2}=(1+\beta_{i})^{2}g_{S^{n-1}}$.
Since locally $g$ is a warped product, we have 
\[
g=dr^{2}+\tilde{f}(r)^{2}g_{H}
\]
near $p_{i}$. By the choice of $p_i$, we may assume that $r(p_i)=0$, and some points go to $p_i$ by the flow of $\partial_r$ or $-\partial_r$. Since the flow of $\partial_r$ gives geodesics, $r$ gives the distance to $p_{i}$ for point near $p_i$. Therefore, for $\rho$ small, 
$dr=d\rho$, and $\tilde{f}^{2}(r)g_{H}=h^{2}(r,\theta)g_{S^{n-1}}$.
Now, 
\[
\left(\lim_{r\to0}\frac{1}{r^{2}}\tilde{f}^{2}(r)\right)g_{H}=\lim_{r\to0}\frac{1}{r^{2}}h^{2}(r,\theta)g_{S^{n-1}}=(1+\beta_{i})^{2}g_{S^{n-1}}.
\]
Thus, $g_{H}=(1+\beta_{i})^{2}c^{2}g_{S^{n-1}}$ for some constant
$c>0$, and we have $g=dr^{2}+f^{2}(r)g_{S^{n-1}}$ for $f(r)=\frac{(1+\beta_{i})}{c}\tilde{f}(r)$.
Moreover, $f'(r)\not=0$ for $r\in(0,\epsilon)$, since $\lim_{r\to0^{+}}f'(r)=1+\beta_{i}$. 

If $p_{i}$ is a regular point of $M$, then it is a critical point
of $u$. On a maximally extended interval $(a,b)$ of $\gamma(r)$,
\[
g=dr^{2}+f(r)^{2}g_{H}.
\]
By the concavity of $f$, the critical point is non-degenerate, hence
they are isolated. In summary, we show that all possible choices of
$p_{i}$ are isolated. Hence, there exist finitely many $p_{i}$. 

Next, we consider the general geodesic flow in $\frac{\partial}{\partial r}$
direction. Since $dr$ is closed, in a neighborhood of $H$, $dr$
is integrable. Adjusting by a constant, we may assume that $r|_{H}=0$.
Then for $\epsilon$ small, $r^{-1}(-\epsilon,\epsilon)$ has a warped
product structure. There exists a $a>0$  such that $\lim_{r\to a^{-}}f(r)=0$.
Since all critical points are isolated, we know that for geodesic
flow $\gamma(x,r)$: 
\[
\lim_{r\to a^{-}}\gamma(x,r)=p,\forall x\in H,
\]
which means that all nearby points flow to one particular choice of
$p_{i}$. 

Similarly, there exists some $b>0$, such that $\lim_{r\rightarrow b^{+}}f(r)=0$.
Note that we may have at most two critical points for warped product
metric.

We have the diffeomorphism 
\[
F:(b,a)\times H\rightarrow M\backslash\{p,q\}
\]
and on $M\backslash\{p,q\},$ by the above argument 
\[
g=dr^{2}+f(r)^{2}g_{H}.
\]

As $g$ is conic metric at $p$, $f(r)^{2}g_{H}=f(r)^{2}g_{S^{n-1}}$
with $\lim_{r\rightarrow0}f(r)/r=1+\beta$, here $\beta$ is the conic
coefficient for $p$.

We summarize our result. For a conic manifold, there exists at least
one conic point and we have three cases as below.

Case A: only one conic point $q$ and without loss of generality $F:(0,+\infty)\times H\rightarrow M\backslash\{q\}$. 

We have
\[
g=dr^{2}+f(r)^{2}g_{S^{n-1}}
\]
 on $M\backslash q$ with $\lim_{r\rightarrow0}f(r)/r=1+\beta$, here
$\beta$ is the conic coefficient of point $q$.

Case B: two critical points $p,q$ including only one conic point
$p$ and $F:(b,a)\times H\rightarrow M\backslash\{p,q\}.$

We have that $g=dr^{2}+f(r)^{2}g_{S^{n-1}}$ on $M\backslash\{q,p\}$
with $\lim_{r\to b}|\frac{f(r)}{r-b}|=1+\beta_{b}$, where $\beta_{b}$
is the the conic coefficient of point $p$, and $\lim_{r\rightarrow a}|\frac{f(r)}{r-a}|=1$.
In this case, due to the existence of regular critical point $p$,
we see that $H$ is topological $S^{n-1}.$

Case C: two conic points $p,q$ and $F:(b,a)\times H\rightarrow M\backslash\{p,q\}.$

We have that $g=dr^{2}+f(r)^{2}g_{S^{n-1}}$ on $M\backslash\{q,p\}$
with $\lim_{r\to b}|\frac{f(r)}{r-b}|=1+\beta_{b}$, where $\beta_{b}$
is the conic coefficient of point $p$, and $\lim_{r\rightarrow a}|\frac{f(r)}{r-a}|=1+\beta_{a},$
where $\beta_{a}$ is the conic coefficient of point $q$. 

In particular, if $g$ is smooth, then $(M,g)$ is conformal to warped
product $N\times(-\infty,+\infty)$, or Euclidean space or a sphere. 
\end{proof}

\section{Conic Andrews inequality }

In this section, we deal with conic manifolds. Let $(M,g_{0},D,g)$
be a conic $n$-manifold with $D=\sum_{i=1}^{k}p_{i}\beta_{i}$. Consider
the following differential equation.
\[
\Delta_{g}u=\phi,
\]
for $\phi\in C^{\alpha}(M,g_{0})$. By the work of F. Wan \cite{wan},
we have a $W^{1,2}(M,g_{0})\cap C^{\alpha'}(M,g_{0})$ weak solution
$u$ for some $\alpha'>0$. By Theorem \ref{thm:main result in appendix}
in the Appendix, we have $u\in C^{1,\gamma_{i}}(B_{\delta}(p_{i}),g)\cap C^{2,\alpha}(M\backslash D,g_{0})$
where $\gamma_{i}$ is given in (\ref{eq:holder index in appendix})
in Theorem \ref{thm:main result in appendix}. Without loss of generality,
we may assume $\gamma_{1}=\min\{\gamma_{i}\}$. Then $u\in C^{1,\gamma_{1}}(M,g)$.

Now, we may give the proof of Theorem \ref{thm:(Conic-Andrews)}.
\begin{proof}
[Proof of Theorem \ref{thm:(Conic-Andrews)}]Let $\nabla$ be the
connection of $g$. For $\int_{M^{n}}\phi dv_{g}=0$ and $\phi\in C^{\alpha}(M,g_{0})$
for some $\gamma>0$, there exists a solution $u\in C^{1,\gamma_{1}}(B_{\delta}(p_{i}),g)\cap C^{\infty}(M\backslash D,g)$
such that $\triangle_{g}u=\phi.$ Let $\varphi_{\varepsilon}$ be
a smooth function such that $supp(1-\varphi_{\varepsilon})\subset\cup_{i=1}^{q}B_{2\varepsilon}(p_{i})$
and $\varphi_{\varepsilon}=0$ on $B_{\varepsilon}(p_{i})$ and $|\partial_{i}\varphi_{\varepsilon}|\le C\frac{1}{\varepsilon}.$

\begin{align*}
 & \lim_{\varepsilon\rightarrow0}\bigg|\int_{M^{n}}\nabla_{i}\varphi_{\varepsilon}\nabla_{j}u\cdot u^{ij}dv_{g}\bigg|\\
 & =\lim_{\varepsilon\rightarrow0}\bigg|\int_{B_{2\varepsilon}\backslash B_{\varepsilon}(p_{i})}\nabla_{i}\varphi_{\varepsilon}\nabla_{j}u\nabla^{ij}udv_{g}\bigg|\\
 & \le C\lim_{\varepsilon\rightarrow0}\bigg|\int_{\varepsilon}^{2\varepsilon}\frac{1}{\varepsilon}\rho^{\gamma_{1}}\rho^{-1+\gamma_{1}}\rho^{n-1}d\rho\bigg|\\
 & =0.
\end{align*}
Similar argument for $\lim_{\varepsilon\rightarrow0}\int_{M^{n}}bu^{k}(\varphi_{\varepsilon})_{k}\triangle uds_{g}=0.$

\begin{align*}
 & \int_{M^{n}\backslash D}|\nabla^{2}u-\frac{\triangle u}{n}g|_{g}^{2}dv_{g}\\
 & =\lim_{\varepsilon\rightarrow0}\int_{M^{n}}\varphi_{\varepsilon}|\nabla^{2}u-\frac{\triangle u}{n}g|_{g}^{2}dv_{g}\\
 & =\lim_{\varepsilon\rightarrow0}\int_{M^{n}}\varphi_{\varepsilon}(\nabla_{i}\nabla_{j}u\nabla^{i}\nabla^{j}u-\frac{(\triangle u)^{2}}{n})dv_{g}\\
 & =\lim_{\varepsilon\rightarrow0}\int_{M^{n}}\varphi_{\varepsilon}\big(\frac{n-1}{n}(\triangle u)^{2}-Ric_{ij}u^{i}u^{j}\big)dv_{g}+\lim_{\varepsilon\rightarrow0}\int_{M^{n}}\nabla_{i}\varphi_{\varepsilon}\nabla_{j}u\cdot u^{ij}dv_{g}\\
 & =\int_{M^{n}\backslash D}\big(\frac{n-1}{n}(\triangle u)^{2}-Ric_{ij}u^{i}u^{j}\big)dv_{g}
\end{align*}
Let $b$ be a constant to be fixed.

\begin{align*}
 & \lim_{\varepsilon\rightarrow0}\int_{M^{n}}\varphi_{\varepsilon}(Ric^{-1})^{ij}(\nabla_{i}\phi+bRic_{ik}u^{k})(\nabla_{j}\phi+bRic_{jl}u^{l})dv_{g}\\
 & =\lim_{\varepsilon\rightarrow0}\int_{M^{n}}\varphi_{\varepsilon}(Ric^{-1})^{ij}\phi_{i}\phi_{j}dv_{g}+2\int_{M^{n}}\varphi_{\varepsilon}bu^{k}\phi_{k}dv_{g}+\int_{M^{n}}\varphi_{\varepsilon}b^{2}Ric_{ij}u^{i}u^{j}dv_{g}\\
 & =\lim_{\varepsilon\rightarrow0}\int_{M^{n}}\varphi_{\varepsilon}(Ric^{-1})^{ij}\phi_{i}\phi_{j}dv_{g}+\int_{M^{n}}\varphi_{\varepsilon}b^{2}Ric_{ij}u^{i}u^{j}v_{g}\\
 & -2\lim_{\varepsilon\rightarrow0}\int_{M^{n}}\varphi_{\varepsilon}b\triangle u\phi dv_{g}-2\lim_{\varepsilon\rightarrow0}\int_{M^{n}}bu^{k}(\varphi_{\varepsilon})_{k}\triangle uds_{g}\\
 & =\lim_{\varepsilon\rightarrow0}\int_{M^{n}}\varphi_{\varepsilon}(Ric^{-1})^{ij}\phi_{i}\phi_{j}dv_{g}+\int_{M^{n}}\varphi_{\varepsilon}b^{2}Ric_{ij}u^{i}u^{j}v_{g}-2\int_{M^{n}}b(\triangle u)^{2}dv_{g}-2\int_{M^{n}}\varphi_{\varepsilon}b\triangle u\phi dv_{g}\\
 & =\int_{M^{n}}(Ric^{-1})^{ij}\phi_{i}\phi_{j}dv_{g}+\int_{M^{n}}b^{2}Ric_{ij}u^{i}u^{j}v_{g}-2\int_{M^{n}}b(\triangle u)^{2}dv_{g}
\end{align*}
Taking $b=1/a$, we obtain

\begin{align*}
0\le\lim_{\varepsilon\rightarrow0} & \int_{M^{n}}\varphi_{\varepsilon}|\nabla^{2}u-\frac{\triangle u}{n}g|_{g}^{2}dv_{g}+a^{2}\int_{M^{n}\backslash D}\varphi_{\varepsilon}(Ric^{-1})^{ij}(\nabla_{i}\phi+bRic_{ik}u^{k})(\nabla_{j}\phi+bRic_{jl}u^{l})\\
= & -\int_{M^{n}\backslash D}Ric_{ij}u^{i}u^{j}dv_{g}+\frac{n-1}{n}\int_{M^{n}\backslash D}(\triangle u)^{2}dv_{g}+a^{2}(\int_{M^{n}\backslash D}(Ric^{-1})^{ij}\phi_{i}\phi_{j}dv_{g}\\
 & +\int_{M^{n}\backslash D}b^{2}Ric_{ij}u^{i}u^{j}dv_{g}-2\int_{M^{n}\backslash D}b\phi^{2}dv_{g})\\
= & \frac{n-1}{n}\int_{M^{n}\backslash D}(\triangle u)^{2}dv_{g}-2a\int_{M^{n}\backslash D}\phi^{2}dv_{g}+a^{2}\int_{M^{n}\backslash D}(Ric^{-1})^{ij}\phi_{i}\phi_{j}dv_{g}.
\end{align*}

So 
\[
(-\frac{n-1}{n}+2a)\frac{1}{a^{2}}\int_{M^{n}\backslash D}\phi^{2}dv_{g}\le\int_{M^{n}\backslash D}(Ric^{-1})^{ij}\phi_{i}\phi_{j}dv_{g}.
\]
Taking $a=\frac{n-1}{n}$, we get the inequality.

When the Andrews inequality (\ref{eq:conic}) is an equality, we have
for $\triangle u=\phi,$

\[
\nabla^{2}u-\frac{\triangle u}{n}g=0\quad on\quad M\backslash D
\]
 and 
\begin{equation}
\nabla_{i}\phi+\frac{n}{n-1}Ric_{ik}u^{k}=0\quad on\quad M\backslash D.\label{eq:second identity}
\end{equation}

\textcolor{black}{By Theorem \ref{thm:warped product-1}, the metric
$g$ can be written as
\begin{equation}
g=f(r)^{2}g_{S^{n-1}}+dr^{2}.\label{eq:add2}
\end{equation}
where $f=u'.$ Since $\int_{M}\triangle udv_{g}=\int_{M}\phi dv_{g}=0,$
we have $f(0)=f(a).$ From Theorem \ref{thm:warped product-1}, we
have $f(0)=f(a)=0$ and have proved this theorem. A direct computation
from (\ref{eq:add2}) shows that the metric is conformally flat.}
\end{proof}
In \cite{Gursky-Streets1}, Gursky-Streets has successfully applied Andrews inequality to study
$\sigma_{2}$-Yamabe problem on smooth four-manifolds. 
It is well known that 
\[
\sigma_{k}(g^{-1}A_{g})=\sum_{1\le i_{1}\cdots<i_{k}\le n}\lambda_{i_{1}}\lambda_{i_{2}}\cdots\lambda_{i_{k}},
\]
where $\{\lambda(A)\}_{i=1}^{n}$ be the set of eigenvalues of $A$
with respect to $g$ and $A=\frac{1}{n-2}(\text{{\rm Ric}}-\frac{{\rm Tr}({\rm Ric})}{2(n-1)}g)$. 

Our rigidity result of Andrews inequality can be further refined if $M$ has constant $\sigma_{\frac{n}{2}}$-curvature. We mention a particular type of conic manifolds that are first discussed
in Chang-Han-Yang \cite{CHY}.
\begin{defn}
$(S^{n},g_{S^{n}},D,g)$ is called a conic $n$-sphere of football
type if $D=\beta p+\beta q$ where $\beta\in(-1,0)$ and $p,q$ are
two distinct points on $S^{n}$.
\end{defn}

\begin{rem}
For a conic $n$-sphere of football type, by a simple conformal transform,
we may assume that $p,q$ are antipodal points on $S^{n}$. It is
also clear that a conic sphere of football type is locally conformally
flat. 
\end{rem}

When $M$ has constant $\sigma_{\frac{n}{2}}$-curvature and the equality in (\ref{eq:conic}) holds, we can prove that $M$ is conic football. For reader's convenience, we give a detailed proof.
\begin{cor}
Under the same conditions for manifolds in Theorem \ref{thm:(Conic-Andrews)},
assuming that $\sigma_{\frac{n}{2}}(g^{-1}A_{g})$ is a constant and
the equality holds in Theorem \ref{thm:(Conic-Andrews)}, we have
$(M,g_{0},D,g)$ is conformal to a conic sphere of football type with
an $O(n-1)$.
\end{cor}

\begin{proof}
By the proof of Theorem \ref{thm:(Conic-Andrews)}, we know
\[
g=f(r)^{2}g_{S^{n-1}}+dr^{2},\quad f(0)=f(a)=0,\quad r\in[0,a].
\]
 By Lemma \ref{lem:petersen}, in polar coordinates, we have 
\end{proof}
\begin{align*}
A_{rr} & =-\frac{f''}{f}-\frac{1}{2}\frac{1-(f')^{2}}{f^{2}}
\end{align*}

\[
A_{ii}=\frac{1}{2}\frac{1-(f')^{2}}{f^{2}},\quad A_{ir}=0.
\]

Therefore by $k=\frac{n}{2}$,

\begin{align*}
\sigma_{k}(A_{g}) & =C_{n-1}^{k}(\frac{1}{2}\frac{1-(f')^{2}}{f^{2}})^{k}+C_{n-1}^{k-1}(\frac{1}{2}\frac{1-(f')^{2}}{f^{2}})^{k-1}(-\frac{f''}{f}-\frac{1}{2}\frac{1-(f')^{2}}{f^{2}})\\
 & =C_{n-1}^{k-1}(\frac{1}{2}\frac{1-(f')^{2}}{f^{2}})^{k-1}(\frac{n-k}{k}(\frac{1}{2}\frac{1-(f')^{2}}{f^{2}})+(-\frac{f''}{f}-\frac{1}{2}\frac{1-(f')^{2}}{f^{2}}))\\
 & =C_{n-1}^{k-1}(\frac{1}{2}\frac{1-(f')^{2}}{f^{2}})^{k-1}(-\frac{f''}{f}).
\end{align*}

Assume $\sigma_{k}(A_{g})=c_{1}.$ Let $c_{0}=c_{1}/(C_{n-1}^{k-1}(\frac{1}{2})^{k-1}).$
Then we have 
\[
\bigg(\frac{1-(f')^{2}}{f^{2}}\bigg)^{k-1}(-\frac{f''}{f})=c_{0}.
\]

Multiplying $f'$ with the equation, there exists a constant $c_{2}$
such that 

\[
\frac{1}{n}(1-(f')^{2})^{k}-\frac{c_{0}}{n}f^{n}=c_{2}.
\]

By $f(0)=f(a)=0$, we have $|f'|(0)=|f'|(a)$. From Theorem \ref{thm:warped product-1},
we know $D=\beta p+\beta q$ where $M\backslash\{p,q\}=S^{n-1}\times(0,a)$.

\section{Andrews inequality for manifolds with boundary}

In this section, we prove Theorem \ref{thm:andrews ineq with Boundary }.
We first prove inequality (\ref{eq:with boundary}).
\begin{proof}
[Proof of Inequality (\ref{eq:with boundary})] Given $\phi\in C^{\infty}(M)$,
we consider the following equation with Neumann boundary condition:
\begin{equation}
\begin{cases}
\Delta u=\phi, & x\in M,\\
u_{\nu}=0, & x\in\partial M.
\end{cases}\label{eq:add3}
\end{equation}
Here (\ref{eq:add3}) is clearly solvable when $\int_{M}\phi=0$.
As in the conic case, we consider the following inequality
\begin{equation}
0\leq\int_{M}\left|\nabla^{2}u-\frac{\Delta u}{n}g\right|^{2}dv_{g}+a^{2}\int_{M}(Ric^{-1})^{ij}(\nabla_{i}\phi+\frac{1}{a}Ric_{i}^{k}u_{k})(\nabla_{j}\phi+\frac{1}{a}Ric_{j}^{k}u_{k})dv_{g}.\label{eq:boundary andrews eq1}
\end{equation}
Note
\begin{equation}
\int_{M}\left|\nabla^{2}u-\frac{\Delta u}{n}g\right|^{2}dv_{g}=\int_{M}|\nabla^{2}u|^{2}dv_{g}-\frac{1}{n}\int_{M}|\Delta u|^{2}dv_{g}.\label{eq:add4}
\end{equation}
By Bochner formula:
\begin{equation}
|\nabla^{2}u|^{2}=\Delta\left(\frac{1}{2}|\nabla u|^{2}\right)-\nabla u\cdot\nabla\Delta u-Ric(\nabla u,\nabla u).\label{eq:add5}
\end{equation}
Combining (\ref{eq:add4}) and (\ref{eq:add5}), we have
\begin{align}
\int_{M}|\nabla^{2}u|^{2} & =\int_{\partial M}\left(\frac{1}{2}\partial_{\nu}|\nabla u|^{2}-\Delta u\partial_{\nu}u\right)d\sigma+\int_{M}|\Delta u|^{2}dv_{g}-\int_{M}Ric(\nabla u,\nabla u)dv_{g},\label{eq:boundary andrews eq1.1}
\end{align}
where $d\sigma$ is the volume form on the boundary, and $\nu$ is
the unit outer normal vector. In the second term in (\ref{eq:boundary andrews eq1}),
we obtain
\begin{align}
a^{2} & \int_{M}(Ric^{-1})^{ij}(\nabla_{i}\phi+\frac{1}{a}Ric_{i}^{k}u_{k})(\nabla_{j}\phi+\frac{1}{a}Ric_{j}^{k}u_{k})dv_{g}\nonumber \\
 & =a^{2}\int_{M}Ric^{-1}(\nabla\phi,\nabla\phi)dv_{g}+2a\int_{M}\nabla\phi\cdot\nabla udv_{g}+\int_{M}Ric(\nabla u,\nabla u)dv_{g}.\label{eq:boundary andrews eq1.2}
\end{align}
Note

\begin{equation}
\int_{M}\nabla\phi\cdot\nabla udv_{g}=\int_{\partial M}\phi\partial_{\nu}ud\sigma-\int_{M}\phi^{2}dv_{g}=\int_{\partial M}\Delta u\partial_{\nu}ud\sigma-\int_{M}\phi^{2}dv_{g}.\label{eq:boundary andrews eq1.3}
\end{equation}
Combining (\ref{eq:boundary andrews eq1}),(\ref{eq:boundary andrews eq1.1}),(\ref{eq:boundary andrews eq1.2}),
and (\ref{eq:boundary andrews eq1.3}), we get
\begin{align}
0 & \leq\left(\frac{n-1}{n}-2a\right)\int_{M}\phi^{2}+a^{2}\int_{M}(Ric^{-1})^{ij}\nabla_{i}\phi\nabla_{j}\phi\nonumber \\
 & \ +\int_{\partial M}\left(\frac{1}{2}\partial_{\nu}|\nabla u|^{2}+(2a-1)\Delta u\partial_{\nu}u\right)d\sigma.\label{eq:boundary andrews ineq 2}
\end{align}
Next, we compute boundary terms. With $u_{\nu}=0$ on $\partial M$,
we write
\begin{equation}
\nabla u|_{\partial M}=u_{i}e_{i},\label{eq:equality}
\end{equation}
where $e_{1},e_{2}\cdots e_{n-1}\in TM$ form a local orthonormal
frame tangential to the boundary. Then 
\begin{align}
\nu(|\nabla u|^{2}) & =2\langle\nabla_{\nu}\nabla u,\nabla u\rangle\label{eq:=00005Cnu(|=00005Cnablau|^2)=00003DHess}\\
 & =2\nabla^{2}u(\nabla u,\nu).\nonumber 
\end{align}
Recalling 
\[
{\rm II}(\alpha,\beta)=\langle\nabla_{\alpha}\beta,-\nu\rangle|_{\partial M}
\]
is the second fundamental form of the boundary $\partial M$ with
respect to the inner normal vector $-\nu$. Since $u_{\nu}=0$, we
obtain
\begin{align}
0 & =e_{k}\langle\nabla u,\nu\rangle=\langle\nabla_{e_{k}}\nabla u,\nu\rangle+\langle\nabla u,\nabla_{e_{k}}\nu\rangle\label{eq:Hess and II connecting equation 1 , sec 4}\\
 & =\nabla^{2}u(e_{i},\nu)+{\rm II}(e_{k},\nabla u).\nonumber 
\end{align}
Thus, by (\ref{eq:equality}) and (\ref{eq:Hess and II connecting equation 1 , sec 4})
\begin{align*}
2\nabla^{2}u(\nabla u,\nu) & =2u_{i}\nabla^{2}u(e_{i},\nu)\\
 & =-2u_{i}{\rm II}(e_{i},\nabla u)\\
 & =-2{\rm II}(\nabla u,\nabla u).
\end{align*}
Notice that the boundary term in (\ref{eq:boundary andrews ineq 2})
is
\begin{align}
\int_{\partial M}\left(\frac{1}{2}\partial_{\nu}|\nabla u|^{2}+(2a-1)\Delta u\partial_{\nu}u\right)d\sigma & =-\int_{\partial M}{\rm II}(\nabla u,\nabla u)d\sigma,\label{eq:coundary term rewrite}
\end{align}
since $u_{\nu}=0$. By (\ref{eq:boundary andrews ineq 2}) and the
assumption ${\rm II}\geq0$, we have 
\[
0\leq\int_{\partial M}{\rm II}(\nabla u,\nabla u)\leq\left(\frac{n-1}{n}-2a\right)\int_{M}\phi^{2}+a^{2}\int_{M}(Ric^{-1})^{ij}\nabla_{i}\phi\nabla_{j}\phi,
\]
Therefore,
\[
-\left(\frac{n-1}{n}-2a\right)\int_{M}\phi^{2}\leq a^{2}\int_{M}(Ric^{-1})^{ij}\nabla_{i}\phi\nabla_{j}\phi.
\]
Then the boundary Andrews inequality (\ref{eq:with boundary}) follows
by choosing $a=\frac{n-1}{n}$.
\end{proof}
Next, we discuss the equality case of (\ref{eq:with boundary}). An
immediate consequence of the proof of (\ref{eq:with boundary}) is
the following corollary:
\begin{cor}
The equality in Andrews inequality (\ref{eq:with boundary}) holds
for a non-constant function $\phi$ if and only if there is a function
$u$ such that 
\begin{equation}
\begin{cases}
\nabla^{2}u-\frac{\Delta u}{n}g & =0,\\
{\rm II}(\nabla u,\nabla u)|_{\partial M} & =0,
\end{cases}\label{eq:equality holds condition}
\end{equation}
where $\Delta u=\phi$. 
\end{cor}

\begin{proof}
By the proof of boundary Andrews inequality, the equality in (\ref{eq:with boundary})
requires that the right hand side of (\ref{eq:boundary andrews eq1})
vanishes for $a=\frac{n-1}{n}$ and ${\rm II}(\nabla u,\nabla u)|_{\partial M}=0$.
Thus, we only need to prove $\nabla\phi+\frac{n}{n-1}Ric(\nabla u)=0$
if (\ref{eq:equality holds condition}) holds. The first restriction
$\nabla^{2}u-\frac{\Delta u}{n}g=0$ implies that $M$ is locally
a warped product. In fact, locally we have the warped product structure,
\[
g_{M}=dr^{2}+f(r)^{2}g_{N},
\]
where $dr=\frac{du}{|du|}$, and $f=|du|$, see Proposition \ref{prop:conformal vector field gives warped product}.
Then, $\phi=f'(r)$ , and by (\ref{eq:Ricc}), 
\[
\nabla\phi+\frac{n}{n-1}Ric(\nabla u)=0.
\]
\end{proof}
Since $u$ restricted on $\partial M$ also yields a locally warped
product structure, we now state a rigidity result for $\partial M$.
\begin{cor}
If the equality in Andrews inequality (\ref{eq:with boundary}) holds
for a non-constant function $\phi$, then $(\partial M,g|_{\partial M})$
is conformal to a $(n-1)$-sphere with the standard structure.\label{cor:boundary is a S_n-1}
\end{cor}

\begin{proof}
Since $\phi\not\equiv0$, $u$ is not a constant. In fact, by Proposition
\ref{prop:critical points have trivial index}, $u$ can not be locally
constant in $M\backslash\partial M$. Let $\bar{u}=u|_{\partial M}$
be the restriction of $u$ on the boundary. Let $\bar{g},\bar{\nabla}$
be the induced metric and connection on $\partial M$. Since $u_{\nu}=0$,
we have that $\overline{\nabla}\bar{u}=\nabla u|_{\partial M}$. Let
$e_{1}\cdots e_{n-1}$ be a local orthonormal frame on $\partial M$.
We compute the Hessian bilinear form of $\bar{u}$:
\begin{align*}
\overline{\nabla}^{2}\bar{u}(e_{i},e_{j}) & =\langle\bar{\nabla}_{e_{i}}\nabla\bar{u},e_{j}\rangle_{\bar{g}}\\
 & =\langle\nabla_{e_{i}}\nabla u-\langle\nabla_{e_{i}}\nabla u,\nu\rangle\nu,e_{j}\rangle_{g}\\
 & =\langle\nabla_{e_{i}}\nabla u,e_{j}\rangle_{g}\\
 & =\nabla^{2}u(e_{i},e_{j})|_{\partial M}.
\end{align*}
Therefore, by (\ref{eq:equality holds condition}),
\[
\overline{\nabla}^{2}\overline{u}=\frac{\Delta u}{n}g|_{\partial M}=\frac{\Delta u}{n}\bar{g}=\frac{\bar{\Delta}\bar{u}}{n-1}\bar{g}.
\]
We claim that $\bar{u}$ is not a constant on each component of $\partial M$.
We argue by contradiction. Let $A$ be a connected component of $\partial M$
and $q\in A$ be a smooth boundary point. Since $u$ is not a constant
and $A$ is a level-set of $u$, Proposition \ref{prop:conformal vector field gives warped product}
shows that in a neighborhood $U$ of $q$, the metric $g$ is locally
a warped product metric $g=dr^{2}+f^{2}(r)g_{A}$ such that $r(q)=0$
and $dr^{2}|_{A\cap U}=0$. Here $g_{A}=g|_{TA}$ is the induced metric
on $A$. If $\bar{u}\equiv c$ and $u_{\nu}=0$ on $A$, then $\nabla u=0$
on $A$ and $f(0)=|\nabla u(q)|=0$. Thus the metric $g_{A}=g|_{TA}=0$
is degenerate in $U\cap A$, which contradicts the fact that $q$
is a smooth boundary point. Therefore, $\bar{u}$ is not a constant
on $\partial M$. We can apply Corollary \ref{cor:warped sphere}
to $\partial M$ and $\bar{u}$ which concludes that each component
of $\partial M$ is a warped sphere.
\end{proof}
If $\partial M$ is connected, by Theorem \ref{thm:original}, $\partial M$
is conformal to the standard $(n-1)$-sphere. In order to finish our
argument for $M$, we need to establish the connect-ness of $\partial M$.
\begin{prop}
\label{prop:partial M is connected}If $M$ is connected, then $u$
can only have one maximum and one minimum. Furthermore, $\partial M$
is connected and $M$ is a hemisphere. 
\end{prop}

\begin{proof}
This is a straightforward application of Morse theory. By Proposition
\ref{prop:critical points have trivial index}, there are only maximums
and minimums on $M$. Suppose that $\partial M$ has at least two
connected components. By the argument in Corollary \ref{cor:boundary is a S_n-1},
$\bar{u}=u|_{\partial M}$ is not a constant on each component of
$\partial M$. Thus, $\bar{u}$ has at least one maximum on each boundary
component. Let $p$ be a local maximum of $\bar{u}$ on $\partial M$,
then $\nabla u|_{\partial M}=\bar{\nabla}\bar{u}=0$, and $u_{\nu}=0$.
Hence, $p$ is a critical point of $u$ in $M$. By Proposition \ref{prop:critical points have trivial index}\textbf{,
$p$ }is a local maximum of $u$ in $M$. Let $p_{1},\cdots,p_{k}$
be maximum points of $u$. Consider the flow given by 
\[
\frac{d}{dt}\phi_{t}(x)=-\nabla u(\phi_{t}(x)).
\]
The flow $\phi$ is compatible with the boundary since $\nabla u$
is tangent to the boundary. Let 
\[
\mathcal{D}(p_{i})=\{x\in M:\lim_{t\to-\infty}\phi_{t}(x)=p_{i}\}.
\]
be the \textbf{unstable manifold} of $p_{i}$. For two different local
maximums $p_{1}$ and $p_{2}$, if critical point $q\in\overline{\mathcal{D}(p_{1})}\cap\overline{\mathcal{D}(p_{2})}$,
then $q$ has to be a critical point of $u$ with non-trivial index
which can not exist by Proposition \ref{prop:critical points have trivial index}.
Therefore, $\overline{\mathcal{D}(p_{1})}\cap\overline{\mathcal{D}(p_{2})}=\emptyset$,
which implies that $M$ is not connected. This is a contradiction.
We have thus proved that $\partial M$ is connected.
\end{proof}
\begin{proof}
[Proof of Theorem \ref{thm:andrews ineq with Boundary }] We prove
the rigidity result when the equality holds in (\ref{eq:with boundary}).
By Proposition \ref{prop:partial M is connected}, there is only one
maximum $p$ and one minimum $q$ on $M$ which both lie on the boundary.
Since $p$ is non-degenerate, near $p$, level-sets of $u$ near local
maximum are diffeomorphic to a $(n-1)$-dimensional hemisphere by
Morse Lemma. Take $\epsilon$ small and let $L$ be the level-set
$u^{-1}(u(p)-\epsilon)$. For a generic choice of $\epsilon$, $L$
is a smooth hypersurface with boundary. Since $\nabla u=\bar{\nabla}u|_{\partial M}$,
$M$ has a warped product metric given by 
\[
g=dr^{2}+f(r)^{2}g_{L},
\]
where $r$ and $f=|\nabla u|$ are globally defined since $M-\{p,q\}$
is contractible. Note $f(r(p))=0$ if $p\in\partial M$ is a critical
point of $u$. Without loss of generality, we assume that $r(p)=0$.
By restricting the warped product structure to $\partial M$, we see
that 
\begin{equation}
\lim_{r\to0}(f(r)^{2}/r^{2})g_{\partial L}=g_{S^{n-2}}.\label{eq:=00005Bf(r)/r=00005D^2g_L}
\end{equation}
In fact, by choosing a geodesic polar coordinate at $p$ on $\partial M$
and $\nabla u=\bar{\nabla}u|_{\partial M}$, we see (\ref{eq:=00005Bf(r)/r=00005D^2g_L})
easily. On the other hand, $p$ is also a maximum point on $M$. Since
$p\in\partial M$, by the same argument, 
\[
\lim_{r\to0}(\frac{f(r)}{r})^{2}g_{L}=g_{S_{+}^{n-1}},
\]
 where $g_{S_{+}^{n-1}}$ is the standard metric on unit hemisphere.
Therefore, $g_{L}=c^{2}\cdot g_{S_{+}^{n-1}}$ for some constant $c$
and the metric on $M$ is given by
\[
g=dr^{2}+\tilde{f}(r)^{2}g_{S_{+}^{n-1}},
\]
for some $\tilde{f}(r)=\frac{1}{c}f(r)$. Now, this warped product
structure can be extended to the whole manifold except the critical
points of $u$. The same argument works for the global minimum point
$q$. Thus, the warped product structure can be extended to $q$.
Finally, a direct computation shows that the second fundamental form
of $\partial M$ vanishes identically, which implies that $\partial M$
is totally geodesic. We have finished the proof.
\end{proof}

\appendix

\section{H${\rm \ddot{o}}$lder regularity of conic Laplacian equation}

In this appendix, we study the H${\rm \ddot{o}}$lder regularity of
conic Laplacian equation. Suppose that $g$ is a conic metric on $M^{n}$
with singularities given by $D=\sum_{i=1}^{k}\beta_{i}p_{i}$, $\beta_{i}\in(-1,0)$.
Consider the following Laplacian equation
\begin{equation}
\Delta_{g}u=\varphi,\label{eq:equation on manifolds}
\end{equation}
$\varphi\in L^{q}(M,g)$ for $q\geq2$. F. Wan \cite{wan} has studied
the solvability of equation (\ref{eq:equation on manifolds}). 
\begin{thm}
[\cite{wan}]\label{thm:wan1 }There exists a solution $u\in W^{1,2}\left(M,g\right)$
satisfying $\triangle_{g}u=\varphi$ for $\varphi\in L^{2}(M,g)$\textup{
if and only if} $\text{ }\int_{M}\varphi dV=0$.
\end{thm}

Near a conic point $p$ with cone index $\beta$, we may write $g=r^{2\beta}g_{0}$
in a neighborhood of $p$ where $r(x)=d_{g_{0}}(x,p)$. Wan proves
the following estimate. 
\begin{thm}
[\cite{wan}]\label{thm:wan2}Suppose $u\in W^{1,2}(M)$ is a solution
of $\triangle_{g}u=\varphi$ and $\varphi\in L^{q}(M,g)$ for some
$q>n/2.$ Then we have for some $s>0$, 
\[
|u(x)-u(y)|\leq C\left(\frac{|x-y|}{r}\right)^{s}\left\{ \left(\frac{1}{r^{n(\text{\ensuremath{\beta}}+1)}}\int_{B_{r}}u^{2}\right)^{\frac{1}{2}}+r^{\left(2-\frac{n}{q}\right)(\beta+1)}\|\varphi\|_{L^{q}\left(B_{r}\right)}\right\} 
\]
and 
\begin{equation}
\sup_{B_{r}}|u|\leq C\left(r^{-\frac{n(\beta+1)}{p}}\|u\|_{L^{p}\left(B_{2r}\right)}+r^{\left(2-\frac{n}{q}\right)(\beta+1)}\|\varphi\|_{L^{q}\left(B_{2r}\right)}\right)\label{eq:C0 estimate}
\end{equation}
for any $x,y\in B_{\frac{r}{2}}$, where $C$ is a positive constant
depending only on $n,q,\beta$ .
\end{thm}

Recall that in Definition \ref{def:conic manifolds}, we have defined
\[
\tilde{C}_{\beta}^{2}(B_{1}(0),g_{0})=\{w\in C^{\infty}(B_{1}(0)\backslash\{0\})\cap C^{0}(B_{1}(0)):|r^{-k\beta}\nabla g_{0}^{k}w|<\infty,k=0,1,2\}.
\]
Let $w\in\tilde{C}_{\beta}^{2}(B_{\delta}(0),g_{0})$ and $g=g_{0}r^{2\beta}e^{2w}$.
In a geodesic polar coordinate, $g_{0}=dr^{2}+f^{2}(r,\theta)d\theta^{2}$,
for some smooth $f=r+O(r^{2})$. By our definition of $w$, we have
\begin{align*}
r^{-\beta}|w_{r}| & <C,r^{-\beta-1}|w_{\theta_{i}}|<C,\\
r^{-2\beta}|w_{rr}|<C, & r^{-2\beta}r^{-1}|w_{r\theta_{i}}|<C,r^{-2\beta-2}|w_{\theta_{i}\theta_{j}}|<C.
\end{align*}

\begin{lem}
\label{lem:geodesic normal coordinate}There is a polar coordinate
in a neighborhood of $0$ such that for some $\delta>0$ and $\rho<\delta$.
\[
g=d\rho^{2}+h(\rho,\theta)d\theta^{2},
\]
where $\|h\|_{C^{2}(B_{\delta}(0)\backslash\{0\})}<\infty$ and $h=(1+\beta)\rho+O\left(\rho^{2}\right)$
. 
\end{lem}

\begin{proof}
Let $\rho(r,\theta)=\int_{0}^{r}t^{\beta}e^{w(t,\theta)}dt$. Then,
$\rho\sim Cr^{1+\beta}$ for small $r$. The Jacobean of $(\rho,\theta_{i})$
is given by
\[
\frac{\partial(\rho,\theta_{i})}{\partial(r,\theta_{i})}=\left(\begin{array}{cc}
r^{\beta}e^{w} & \int_{0}^{r}t^{\beta}e^{w}w_{\theta_{i}}dt\\
0 & {\rm Id}
\end{array}\right),
\]
and the inverse Jacobean is given by
\[
\left(\frac{\partial(\rho,\theta_{i})}{\partial(r,\theta_{i})}\right)^{-1}=\left(\begin{array}{cc}
r^{-\beta}e^{-w} & -r^{-\beta}e^{-w}\int_{0}^{r}t^{\beta}e^{w}w_{\theta_{i}}dt\\
0 & {\rm Id}
\end{array}\right)=\left(\begin{array}{cc}
r^{-\beta}e^{-w} & O(\rho^{2-\frac{\beta}{1+\beta}})\\
0 & {\rm Id}
\end{array}\right).
\]
Thus, for $\rho$ small, $(\rho,\theta_{i})$ gives a coordinate near
the origin. Let $h(\rho,\theta):=e^{w}r^{\beta}f(r,\theta)$. We see
that 
\[
g=d\rho^{2}+h^{2}(\rho,\theta)d\theta^{2}.
\]
Then $h$ is smooth in $B_{\delta'}(0)\backslash\{0\}$ and we find
\begin{align*}
h_{\rho} & =f_{r}+\beta r^{-1}f+fw_{r}+O(\rho^{2-\frac{\beta}{1+\beta}})r^{\beta}\sum_{i=1}^{n-1}\left(w_{\theta_{i}}f+wf_{\theta_{i}}\right)\\
 & =1+\beta+O(\rho)+O(\rho^{3}),\\
h_{\theta_{i}} & =e^{w}r^{\beta}\left(w_{\theta_{i}}f+f_{\theta_{i}}\right)=O(\rho),\\
h_{\rho\rho} & =e^{-w}r^{-\beta}\left(f_{rr}-r^{-2}f\beta+r^{-1}f_{r}\beta+f_{r}w_{r}+fw_{rr}\right)+O(\rho^{2})\\
 & =e^{-w}\left(f_{rr}r^{-\beta}-O(r^{-\beta})+f_{r}w_{r}r^{-\beta}+w_{rr}fr^{-\beta}\right)+O(\rho^{2})\\
 & =e^{-w}\left(O(r^{-\beta})+O(1)+O(\rho)\right)=O(1),\\
h_{\rho\theta_{i}} & =f_{r\theta_{i}}+\beta r^{-1}f_{\theta_{i}}+f_{\theta_{i}}w_{r}+fw_{r\theta_{i}}+O(\rho^{2})=O(\rho),\\
h_{\theta_{i}\theta_{j}} & =e^{w}r^{\beta}\left(w_{\theta_{i}}w_{\theta_{j}}f+w_{\theta_{i}\theta_{j}}f+w_{\theta_{i}}f_{\theta_{j}}+w_{\theta_{j}}f_{\theta_{i}}+f_{\theta_{i}\theta_{j}}\right)\\
 & =O(\rho^{2}).
\end{align*}
In particular, we have $h=(1+\beta)\rho+O(\rho^{2})$. 
\end{proof}
Consider the following equation

\begin{equation}
\triangle_{g}u(x)=\varphi(x),\quad x\in B_{1}(0),\label{eq:equation for u}
\end{equation}
where $\varphi\in C^{\gamma}(M,g_{0})$ for some $\gamma>0$. We first
prove an interior regularity result.
\begin{lem}
\label{lem:Schauder} If $u\in C^{s}(B_{1}(0),g_{0})$ is a solution
of equation (\ref{eq:equation for u}), then there is a positive number
$\alpha=\min\{s,2+2\beta\}$ such that for some $\delta>0$ and $\rho<\delta$
\begin{equation}
|u-u(0)|+\rho|\nabla_{g}u|+\rho^{2}|\nabla_{g}^{2}u|<C\rho^{\frac{\alpha}{1+\beta}}.\label{eq:Schauder estimate}
\end{equation}
Under $g_{0}$, for some $\delta'>0$ and $0<r<\delta'$ , we have
\[
|u-u(0)|+r|\nabla_{g_{0}}u|+r^{2}|\nabla_{g_{0}}^{2}u|<Cr^{\alpha}.
\]
\end{lem}

\begin{proof}
We can write $g=d\rho^{2}+h^{2}(\rho,\theta)d\theta^{2}$, where $h(\rho,\theta)=(1+\beta)\rho+O(\rho^{2})$
and is smooth. Therefore, if $\triangle_{g}u(r,\theta)=\varphi(r,\theta)$,
we have
\[
u_{\rho\rho}+\frac{n-1}{h}h_{\rho}u_{\rho}+\frac{\triangle_{S^{n-1}}u}{h^{2}}+\frac{n-3}{h^{3}}\langle\nabla_{\theta}h,\nabla_{\theta}u\rangle_{S^{n-1}}=\varphi(r(\rho),\theta).
\]
Differing by a term $c_{1}\rho^{2}+c_{2}$, we may assume that $u(0,\theta)=0$
and $\varphi(0,\theta)=0$. Clearly, $u_{0}(\rho,\theta):=c_{1}\rho^{2}+c_{2}$
satisfies the estimate \ref{eq:Schauder estimate}. Let $\rho=e^{t}.$
Then
\[
u_{\rho}=e^{-t}u_{t},u_{\rho\rho}=e^{-2t}u_{tt}-e^{-2t}u_{t}.
\]
So for $(t,\theta)\in(-\infty,0)\times S^{n-1}$, we have
\[
u_{tt}+\frac{e^{2t}}{h^{2}}\triangle_{S^{n-1}}u+\left((n-1)\frac{e^{t}h_{\rho}}{h}-1\right)u_{t}+(n-3)\frac{e^{2t}}{h^{3}}\nabla_{\theta}h\cdot\nabla_{\theta}u=\varphi e^{2t}.
\]
Let $\tilde{\varphi}=\varphi e^{2t}$. Let $\Omega_{k}=\{x\in B_{1}(0):e^{-k-1}<\rho<e^{-k}\}$
and $\tilde{\Omega}_{k}=\Omega_{k+1}\cup\Omega_{k}\cup\Omega_{k-1}$.
It is easy to check that $e^{2t}/h^{2},\frac{e^{t}h_{\rho}}{h},$
and $\frac{e^{2t}}{h^{3}}\nabla_{\theta}h$ are all $C^{1}(B_{1}(0)\backslash\{0\})$
functions and their $C^{1}$ norms do not depend on $t$. Then, by
interior Schauder estimate, we have
\begin{equation}
\|u\|_{C^{2,\alpha'}(\Omega_{k})}\leq C(\|u\|_{C^{0}(\tilde{\Omega}_{k})}+\|\tilde{\varphi}\|_{C^{\alpha'}(\tilde{\Omega}_{k})}),\label{eq:Schauder interior estimate}
\end{equation}
for some $\alpha'<\min\{\gamma,1\}$ and $C=C(n,\beta)$. Since we
have assumed that $u,\varphi=0$ at the origin, $|u(e^{t},\theta)|<Ce^{\frac{s}{1+\beta}t}$,
and $|\varphi(r,\theta)|<C'e^{\frac{\gamma}{1+\beta}t}$. Note that
$\varphi$ is also $C^{\gamma}(B_{1}(0),g)$. Let $x_{1}=(e^{t_{1}},\theta_{1}),x_{2}=(e^{t_{2}},\theta_{2})$
be two points in $\tilde{\Omega}_{k}$. Then
\begin{align*}
|\tilde{\varphi}(x_{1})-\tilde{\varphi}(x_{2})| & =|\varphi(x_{1})(e^{2t_{1}}-e^{2t_{2}})+(\varphi(x_{1})-\varphi(x_{2}))e^{2t_{2}}|\\
 & \leq C(|\varphi(x_{1})|e^{-k}d_{g}(x_{1},x_{2})+d_{g}(x_{1},x_{2})^{\gamma}e^{-2k})\\
 & \leq C\left(e^{-k(1+\frac{\gamma}{1+\beta})}d_{g}(x_{1},x_{2})+e^{-2k}d_{g}(x_{1},x_{2})^{\gamma}\right)\\
 & \leq C\left(e^{-k(2+\frac{\gamma}{1+\beta}-\alpha')}+e^{-k(2+\gamma-\alpha')}\right)d_{g}(x_{1},x_{2})^{\alpha'}
\end{align*}
Hence $\|\tilde{\varphi}\|_{C^{\alpha'}(\tilde{\Omega}_{k})}<Ce^{-k(2+\alpha'')}$
for some $\alpha''<\text{\ensuremath{\gamma-\alpha'}}$. Therefore,
by (\ref{eq:Schauder interior estimate}), we have 
\[
\|u\|_{C^{2,\alpha'}(\Omega_{k})}\leq Ce^{-k\alpha'''},
\]
for some $\alpha'''=\min\{\frac{s}{1+\beta},2\}$. Thus 
\[
|u|+|\nabla_{(t,\theta)}u|+|\nabla_{(t,\theta)}^{2}u|<Ce^{t\alpha'''}.
\]
So if we change back to polar coordinate, we have 
\[
|u|+\rho|\nabla_{g}u|+\rho^{2}|\nabla_{g}^{2}u|<C\rho^{\alpha'''}.
\]
If we change back to $g_{0}$ metric, then, 
\[
|u|+r|\nabla_{g_{0}}u|+r^{2}|\nabla_{g_{0}}^{2}u|<Cr^{\alpha},
\]
for $\alpha=\min\{s,2+2\beta\}$.
\end{proof}
Now, we can state our main result in this appendix. 
\begin{thm}
\label{thm:main result in appendix}Let $u$ be the solution to $\triangle_{g}u=\varphi$.
Near a conic singularity with index $\beta$, for some small $\delta>0$,
$u\in C^{1,\alpha}(B_{\delta}(0),g)$ for
\begin{equation}
\alpha=\begin{cases}
\tau, & \rm{for\ any}\ 0<\tau<\min\{\gamma_{1}-1,1\},\quad{\rm if}\ -\frac{1}{2}<\beta<0,\\
1, & {\rm if}\ -1<\beta\leq-\frac{1}{2},
\end{cases}\label{eq:holder index in appendix}
\end{equation}
where 
\[
\gamma_{1}=\sqrt{(\frac{n-2}{2})^{2}+\frac{n-1}{(1+\beta)^{2}}}-\frac{n-2}{2}>1.
\]
Moreover, under metric $g_{0}$, we have
\begin{equation}
|u-u(0)|+r|\nabla_{g_{0}}u|+r^{2}|\nabla_{g_{0}}^{2}u|<Cr^{(1+\alpha(\beta))(1+\beta)}.\label{eq:Schauder estimate in theorem}
\end{equation}
\end{thm}

\begin{proof}
We write equation (\ref{eq:equation for u}) as 
\begin{equation}
u_{\rho\rho}+\frac{n-1}{h}h_{\rho}u_{\rho}+\frac{\triangle_{S^{n-1}}u}{h^{2}}+\frac{n-3}{h^{3}}\langle\nabla_{\theta}h,\nabla_{\theta}u\rangle_{S^{n-1}}=\varphi.\label{eq:1}
\end{equation}
Let $\triangle_{c}$ be an operator
\[
\triangle_{c}u=u_{\rho\rho}+\frac{n-1}{\rho}u_{\rho}+\frac{\triangle_{S^{n-1}}u}{(1+\beta)^{2}\rho^{2}}.
\]
Then we rewrite (\ref{eq:1}) as 
\begin{equation}
\triangle_{c}u=\varphi+(\triangle_{g}-\triangle_{c})u=\varphi+O(1)u_{\rho}+O(\rho^{-1})(|\nabla_{\theta}u|+|\triangle_{S^{n-1}}u|).\label{eq:=00005CDetla_cone u equation}
\end{equation}
We notice that differing a term in the form of $c_{1}\rho^{2}+c_{2}$,
we may assume that $u(0)$ and $\varphi(0)$ are both equal to zero
without changing the form of equation (\ref{eq:=00005CDetla_cone u equation})
and our regularity result. In particular, under this assumption, we
have $|u|<C\rho^{\frac{s}{1+\beta}}$ and $|\varphi|<C\rho^{\frac{\gamma}{1+\beta}}$.
We assume without loss of generality that $s\leq2+2\beta$.

We first prove that $|\nabla_{\theta}u|\le C\rho^{\frac{s}{1+\beta}}$
and $|\triangle_{S^{n-1}}u|\le C\rho^{\frac{s}{1+\beta}}$ if $u$
is $C^{s}$ H$\ddot{o}$lder continuous. We rewrite (\ref{eq:1})
and by Lemma \ref{lem:Schauder}, 
\begin{align}
 & \triangle_{S^{n-1}}u+\frac{(n-3)}{h}\langle\nabla_{\theta}h,\nabla_{\theta}u\rangle_{S^{n-1}}\label{eq:tangent direction}\\
 & =-h^{2}(u_{\rho\rho}+\frac{(n-1)h_{\rho}}{h}u_{\rho}-\varphi)\nonumber \\
 & =O(\rho^{\frac{s}{1+\beta}})+O(\rho^{2+\frac{\gamma}{1+\beta}}).\nonumber 
\end{align}
Now by $W^{2,p}$ estimates for (\ref{eq:tangent direction}), we
know that $|u|_{C^{1,\gamma_{3}}(S^{n-1})}\le C$ for any $0<\gamma_{3}<1$.
Since $|u|=O(\rho^{\frac{s}{1+\beta}})$, by interpolation inequality
of H${\rm \ddot{o}}$lder norm, we have 
\begin{align*}
|\nabla_{\theta}u|_{C^{0}(S^{n-1})} & \le C|u|_{C^{1,\gamma_{3}}(S^{n-1})}^{\frac{1}{1+\gamma_{3}}}|u|_{C^{0}(S^{n-1})}\\
 & \le O(\rho^{\frac{s}{1+\beta}}),
\end{align*}
and then by (\ref{eq:1})
\[
|\triangle_{S^{n-1}}u|\le O(\rho^{\frac{s}{1+\beta}})+O(\rho^{2+\frac{\gamma}{1+\beta}}).
\]
Let
\begin{equation}
\widetilde{\varphi}:=\varphi+(\triangle_{g}-\triangle_{c})u=O(\rho^{\min\{\frac{s}{1+\beta}-1,\frac{\gamma}{1+\beta}\}}).\label{eq:tilde =00005Cvarphi}
\end{equation}
Then (\ref{eq:1}) becomes

\begin{equation}
\triangle_{c}u=\widetilde{\varphi}.\label{eq:3}
\end{equation}

Next, we consider the spherical harmonic expansion of $u$. Let $u(\rho,\theta)=\sum_{j=0}^{\infty}u_{j}(\rho)\psi_{j}(\theta)$,
where $\triangle_{S^{n-1}}\psi_{j}(\theta)=-\lambda_{j}\psi_{j}(\theta)$
and $\lambda_{j}$ is the $j^{th}$ eigenvalue of $\triangle_{S^{n-1}}$.
We have 

\begin{align}
\triangle_{c}(\sum_{i=0}^{\infty}u_{i}(\rho)\psi_{i}(\theta)) & =\sum_{i=0}^{\infty}\widetilde{\varphi}_{i}(\rho)\psi_{i}(\theta).\label{eq:decompostion}
\end{align}
Then $u_{i}$ satisfies the equation:

\[
u_{i}''(\rho)+\frac{(n-1)}{\rho}u_{i}'(\rho)-\frac{\lambda_{i}}{(1+\beta)^{2}\rho^{2}}u_{i}=\widetilde{\varphi}_{i}.
\]
The corresponding homogenous equation has two solutions $\rho^{\alpha_{i}^{+}},$
$\rho^{\alpha_{i}^{-}}$, where 
\[
\alpha_{i}^{+}=\sqrt{(\frac{n-2}{2})^{2}+\frac{\lambda_{i}}{(1+\beta)^{2}}}-\frac{n-2}{2},
\]

\[
\alpha_{i}^{-}=-\sqrt{(\frac{n-2}{2})^{2}+\frac{\lambda_{i}}{(1+\beta)^{2}}}-\frac{n-2}{2}.
\]
Therefore by the H$\ddot{o}$lder continuity of $u$, 
\begin{align}
u_{i}(\rho) & =c\rho^{\alpha_{i}^{+}}-\rho^{\alpha_{i}^{+}}\int_{1}^{\rho}\frac{t^{\alpha_{i}^{-}}\tilde{\varphi}_{i}}{(\alpha_{i}^{+}-\alpha_{i}^{-})t^{\alpha_{i}^{+}+\alpha_{i}^{-}-1}}dt+\rho^{\alpha_{i}^{-}}\int_{0}^{\rho}\frac{t^{\alpha_{i}^{+}}\tilde{\varphi}_{i}}{(\alpha_{i}^{+}-\alpha_{i}^{-})t^{\alpha_{i}^{+}+\alpha_{i}^{-}-1}}dt\label{eq:for u_i}\\
 & =:c\rho^{\alpha_{i}^{+}}+v_{i}(\rho).\nonumber 
\end{align}
For $\lambda_{0}=0$, $\alpha_{0}^{+}=0$ , $\alpha_{0}^{-}=(2-n)$,
and
\begin{align*}
u_{0}(\rho) & =c-\int_{1}^{\rho}\frac{t^{\alpha_{0}^{-}}\tilde{\varphi}_{i}}{-\alpha_{0}^{-}t^{\alpha_{0}^{-}-1}}dt+\rho^{\alpha_{0}^{-}}\int_{0}^{\rho}\frac{\tilde{\varphi}_{i}}{-\alpha_{0}^{-}t^{\alpha_{0}^{-}-1}}dt\\
 & =c-\int_{1}^{\rho}\frac{t\tilde{\varphi}_{i}}{-\alpha_{0}^{-}}dt+\rho^{\alpha_{0}^{-}}\int_{0}^{\rho}\frac{\tilde{\varphi}_{i}}{-\alpha_{0}^{-}t^{\alpha_{0}^{-}-1}}dt.
\end{align*}
Since $|\tilde{\varphi}_{i}|\le O(\rho^{\frac{s}{1+\beta}-1})+O(\rho^{\frac{\gamma}{1+\beta}})$,
we know $u_{0}\in C^{\gamma_{0}}(B_\delta(0),g)$, where $\gamma_{0}=\min\{\frac{s}{1+\beta}+1,\frac{\gamma}{1+\beta}+2\}$.
Note $u_{0}\in C^{1,\gamma_{0}-1}(B_\delta(0),g)$ if $\gamma_{0}>1$.

Let $\hat{u}:=u-\sum_{i=0}^{n}u_{i}(\rho)\psi_{i}(\theta)$.

\[
\int_{\mathbb{S}^{n-1}}\widehat{u}\triangle_{\theta}ud\theta=-\int_{\mathbb{S}^{n-1}}\left|\nabla_{\theta}\widehat{u}\right|^{2}d\theta,
\]
and 
\[
\int_{\mathbb{S}^{n-1}}u_{\rho}\widehat{u}d\theta=\int_{\mathbb{S}^{n-1}}\widehat{u}_{\rho}\widehat{u}d\theta=\frac{1}{2}\frac{d}{d\rho}\int_{\mathbb{S}^{n-1}}|\widehat{u}|^{2}d\theta.
\]
Denote $\hat{\tilde{\varphi}}=\widetilde{\varphi}-\sum_{i=0}^{n}\widetilde{\varphi}_{i}(\rho)\psi_{i}(\theta)$
and $y(\rho)^{2}=\int_{S^{n-1}}|\hat{u}|^{2}d\theta$, it holds that 

\[
yy^{\prime\prime}=\int_{\mathbb{S}^{n-1}}\left\{ \widehat{u}_{\rho\rho}\widehat{u}+\left|\widehat{u}_{\rho}\right|^{2}\right\} d\theta-\left|y^{\prime}\right|^{2}.
\]
By (\ref{eq:decompostion}), it holds that 

\begin{align*}
yy''+|y'|^{2}-\int_{S^{n-1}}|\widehat{u}_{\rho}|^{2}d\theta+\frac{(n-1)}{2\rho}\frac{d}{d\rho}(\int_{S^{n-1}}\widehat{u}^{2}d\theta)-\frac{\int_{S^{n-1}}|\nabla_{\theta}\widehat{u}|^{2}d\theta}{(1+\beta)^{2}\rho^{2}}\ \\
\ \ =\int_{S^{n-1}}\widehat{\widetilde{\varphi}}\widehat{u}d\theta & ,
\end{align*}
and then 
\begin{align*}
yy''+\frac{(n-1)}{\rho}yy' & \ge\frac{1}{(1+\beta)^{2}\rho^{2}}\int_{S^{n-1}}|\nabla_{\theta}\widehat{u}|^{2}d\theta+\int_{S^{n-1}}\widehat{\widetilde{\varphi}}\widehat{u}d\theta\\
 & \ge\frac{2n}{(1+\beta)^{2}\rho^{2}}\int_{S^{n-1}}|\widehat{u}|^{2}d\theta-\big(\int_{S^{n-1}}|\widehat{\widetilde{\varphi}}|^{2}\big)^{\frac{1}{2}}\left(\int_{S^{n-1}}|\widehat{u}|^{2}d\theta\right)^{\frac{1}{2}}\\
 & \ge\frac{2n}{(1+\beta)^{2}\rho^{2}}y^{2}-\big(\int_{S^{n-1}}|\widehat{\widetilde{\varphi}}|^{2}\big)^{1/2}y.
\end{align*}
Here we have use the fact that $\int_{S^{n-1}}|\nabla_{\theta}\hat{u}|^{2}d\theta\geq2n\int_{S^{n-1}}|\hat{u}|^{2}d\theta$.
We obtain that 

\[
y''+\frac{(n-1)}{\rho}y'-\frac{2n}{(1+\beta)^{2}\rho^{2}}y\ge-\big(\int_{S^{n-1}}|\widehat{\widetilde{\varphi}}|^{2}\big)^{1/2}.
\]
To continue the proof, we consider two cases: 

$\spadesuit$Case 1: $0>2\beta>-1$.

We have $u_{0}\in C^{1,\gamma_{0}-1}(B_\delta(0),g)$ . For $i=1,2,\cdots,n$, we
have $\lambda_{i}=n-1$. Let 
\[
\gamma_{1}=\alpha_{1}^{+}=\sqrt{(\frac{n-2}{2})^{2}+\frac{n-1}{(1+\beta)^{2}}}-\frac{n-2}{2}.
\]
Then we have $\gamma_{1}>\frac{1}{1+\beta}>1$. By (\ref{eq:tilde =00005Cvarphi}),
$|\tilde{\varphi}_{i}|\le O(\rho^{-1+\hat{\alpha}})$ for $\hat{\alpha}=\min\{\frac{s}{1+\beta},\frac{\gamma}{1+\beta}+1\}>0$.
When $\gamma_{1}\not=1+\hat{\alpha},$ then by (\ref{eq:for u_i}),

\[
|v_{i}(\rho)|\le C\rho^{1+\hat{\alpha}}.
\]
 When $\gamma_{1}=1+\hat{\alpha}$, $|v_{i}|\le C\rho^{1+\hat{\alpha}}|\ln\rho|$.
So in general, we have $u_{i}\in C^{1,\gamma_{2}}$ for $1<\gamma_{2}<\min\{\gamma_{1}-1,\hat{\alpha}\}$. 

Now let us check that the regularity of the remaining term $\hat{u}=u-\sum_{i=0}^{n}u_{i}(\rho)\psi_{i}(\theta)$.
For $i\geq n+1$, $\lambda_{i}\geq2n$. Define $z(\rho)$ by the equation:
\[
z''+\frac{(n-1)}{\rho}z'-\frac{2n}{(1+\beta)^{2}\rho^{2}}z=-\big(\int_{S^{n-1}}|\widehat{\widetilde{\varphi}}|^{2}\big)^{1/2}.
\]
Then 

\[
z=C\rho^{\alpha_{i}^{+}}-\rho^{\alpha_{i}^{+}}\int_{1}^{\rho}\frac{t^{\alpha_{i}^{-}}||\widehat{\widetilde{\varphi}}(\rho,\theta)||_{L^{2}(S^{n-1})}}{(\alpha_{i}^{+}-\alpha_{i}^{-})t^{\alpha_{i}^{+}+\alpha_{i}^{-}-1}}dt+\rho^{\alpha_{i}^{-}}\int_{0}^{\rho}\frac{t^{\alpha_{i}^{+}}||\widehat{\widetilde{\varphi}}(\rho,\theta)||_{L^{2}(S^{n-1})}}{(\alpha_{i}^{+}-\alpha_{i}^{-})t^{\alpha_{i}^{+}+\alpha_{i}^{-}-1}}dt
\]
where $\alpha_{i}^{+}>2$ for $i=n+1$, $\lambda_{i}=2n$, and $||\widehat{\widetilde{\varphi}}||_{L^{2}(S^{n-1})}\le C\rho^{\hat{\alpha}-1}$
, which yields that $z\le C\rho^{1+\hat{\alpha}}.$ If $a$ is sufficiently
large, $az$ satisfies 

\[
(az)''+\frac{(n-1)}{\rho}(az)'-\frac{2n}{(1+\beta)^{2}\rho^{2}}az=-a\big(\int_{S^{n-1}}|\widehat{\widetilde{\varphi}}|^{2}\big)^{1/2}<-\big(\int_{S^{n-1}}|\widehat{\widetilde{\varphi}}|^{2}\big)^{1/2}
\]
and 
\[
az(1)>y(1).
\]
Since $az(0)=y(0)$, by the maximum principle, we have for $\rho\in(0,1)$
\[
az\ge y.
\]
Thus, 
\begin{align*}
\|\hat{u}\|_{L^{2}(B_{\rho})}^{2} & \leq C{\rm vol}(S^{n-1})\int_{0}^{\rho}y^{2}(t)t^{n-1}dt\\
 & \leq C\int_{0}^{\rho}t^{n-1+(2+2\hat{\alpha})}dt\\
 & =C\rho^{n+2+2\hat{\alpha}}
\end{align*}
By Wan' s argument in \cite{wan}, for any $\rho$ small, 

\begin{align*}
\sup_{B_{\rho}}|\widehat{u}| & \le C(\rho^{-\frac{n}{2}}||\widehat{u}||_{L^{2}(B_{2\rho})}+\rho^{2}||\widehat{\widetilde{\varphi}}||_{L^{\infty}(B_{2\rho})})\\
 & \le C(\rho^{-\frac{n}{2}}(\int_{0}^{\rho}t^{n-1}t^{2(1+\hat{\alpha})}dt)^{1/2}+\rho^{2}||\widehat{\widetilde{\varphi}}||_{L^{\infty}(B_{2\rho})})\\
 & \le C\rho^{1+\hat{\alpha}}
\end{align*}
Therefore we obtain that $u\in C^{1,\gamma_{3}}(B_{\delta}(0),g)$,
where $0<\gamma_{3}<\min\{\gamma_{1}-1,\hat{\alpha}\}=\min\{\gamma_{1}-1,\frac{s}{1+\beta},\frac{\gamma}{1+\beta}+1\}$.
Now we can improve the regularity of $u.$ Since $u\in C^{1,\gamma_{3}}$,
we have $|\widetilde{\varphi}|\le C\rho^{\min\{\frac{\gamma}{1+\beta},1\}}$
and $|u|<C\rho$. Now we can apply the above argument with $s=1+\beta$
again to obtain $u\in C^{1,\gamma_{4}}$ for any $\gamma_{4}<\min\{\text{\ensuremath{\gamma_{1}}}-1,1\}.$ 

$\spadesuit$ Case 2: $-2<2\beta\leq-1$.

For $i=1,2\cdots,n$, $\lambda_{i}=n-1$, $\alpha_{1}^{+}>2$. Since
$\alpha_{1}^{+}>2$, $|v_{i}(\rho)|\le C\rho^{1+\hat{\alpha}}$ and
then $u_{i}\in C^{1,\hat{\alpha}}(B_{\delta}(0),g),$ for $i=1,2,\cdots n$
and $\delta>0$. Next, we consider for $\hat{u}=u-\sum_{i=0}^{n}u_{i}$.
We can use the similar argument as in Case 1. Define $z$ by

\[
z''+\frac{(n-1)}{\rho}z'-\frac{2n}{(1+\beta)^{2}\rho^{2}}z=-\big(\int_{S^{n-1}}|\widehat{\widetilde{\varphi}}|^{2}\big)^{1/2}.
\]
\[
az\ge y.
\]
 Then we can use the similar argument as in Case 1 to show that $|z|<C\rho^{1+\hat{\alpha}}$,
for $\hat{\alpha}=\min\{\frac{s}{1+\beta},\frac{\gamma}{1+\beta}+1\}$.
By (\ref{eq:C0 estimate}), for small $\rho$, 

\begin{align*}
\sup_{B_{\rho}}|\widehat{u}| & \le C(\rho^{-\frac{n}{2}}||\widehat{u}||_{L^{2}(B_{2\rho})}+\rho^{2}||\widehat{\widetilde{\varphi}}||_{L^{\infty}(B_{2\rho})})\\
 & \le C(\rho^{-\frac{n}{2}}\big(\int_{0}^{\rho}t^{n-1}t^{2(2\beta+2)}dt)^{1/2}+\rho^{2}||\widehat{\widetilde{\varphi}}||_{L^{\infty}(B_{2\rho})})\\
 & \le C\rho^{1+\hat{\alpha}}.
\end{align*}
Recall that we have proved
$u_{0}\in C^{\gamma_{0}}(B_\delta(0),g)$ and $u_{i}\in C^{1,\hat{\alpha}}(B_{\delta}(0),g)$.
Since $\gamma_{0}=1+\hat{\alpha}$ , we have $u\in C^{1,\hat{\alpha}}(B_\delta(0),g)$.
We can then apply $s=1+\beta$ in the above argument again to obtain
that $u\in C^{1,1}(B_{\delta}(0),g)$. 

Now, we can apply Lemma \ref{lem:Schauder} to obtain (\ref{eq:Schauder estimate in theorem}).
We have finished the proof.

\bibliographystyle{plain}
\bibliography{andrews-fmw}
\end{proof}

\end{document}